\documentclass{amsart}

\title{Geometric quasi-isometric embeddings into Thompson's group $F$}

\author{Sean Cleary} 
\author{Jennifer Taback}

\thanks{The first author acknowledges support from PSC-CUNY grant
    \#63438-0032.  The second author acknowledges partial support from
the NSF-AWM Mentoring Travel Grant and would like to thank the University of Utah for
its hospitality during the writing of this paper.}

\keywords{Thompson's group F, quasi-isometric embeddings}

\subjclass{20F32}

\newcommand\Z{\mathbb Z}

\newcommand\F{\mathcal F}

\newcommand\x{x_0}

\newcommand\Ll{L_L}
\newcommand\rni{R_{NI}}

\newif\ifpdf
    \ifx\pdfoutput\undefined
    \pdffalse 
    \else
    \pdfoutput=1 
    \pdftrue
    \fi

    \ifpdf
    \usepackage[pdftex]{graphicx}
    \else
    \usepackage{graphicx}
    \fi

\newtheorem{theorem}{Theorem}[section]

\newtheorem{lemma}[theorem]{Lemma}
\newtheorem{corollary}[theorem]{Corollary}

\begin{document} 

\begin{abstract}  
We use geometric techniques to investigate several examples
of quasi-isometrically embedded subgroups of  Thompson's group $F$. 
Many of these are explored using the metric properties 
of the shift map $\phi$ in $F$.
These subgroups have simple geometric but complicated algebraic descriptions.
We present them to illustrate the intricate geometry of Thompson's group $F$ as well
as the interplay between its standard finite and infinite presentations.  
These subgroups include those of the form $F^m \times \Z^n$, for integral $m,n \geq 0$,
which were shown to occur as quasi-isometrically embedded subgroups by Burillo
and Guba and Sapir.
\end{abstract} 

\maketitle

\section{Introduction}
\label{sec:intro}

Bridson asked the question of whether  a quasi-isometry exists between Thompson's group $F$ and the group $F
\times \Z$.
Burillo \cite{burillo} provides an example of a quasi-isometric embedding of $F \times \Z$
into $F$ as possible evidence addressing this question. We state his
example in \S \ref{sec:examples} below.
While investigating this question, we came across some interesting examples
of quasi-isometric embeddings into $F$ which we describe below.
These quasi-isometric embeddings all have simple geometric interpretations,
which are often easier to express than the corresponding algebraic or group
theoretic definitions.
We present these examples to illustrate the beautiful geometry evident in Thompson's group $F$.  Our examples are based on the interaction between the finite and infinite presentations of $F$ and the representation of elements of $F$ as pairs of binary rooted trees.  These embeddings
use  shift maps of $F$ and provide concrete geometric
realizations of subgroups of $F$ of the form $F^m \times Z^n$ which are known to be quasi-isometricly embedded by work of  Burillo \cite{burillo} and Guba and Sapir \cite{diag,diag2}
\section{Thompson's group $F$}
\label{sec:thompson}

Thompson's group $F$ has both finite and infinite presentations; it is
usually presented finitely as
$$
{\mathcal F} = \langle x_0,x_1 |
[x_0x_1^{-1},x_0^{-1}x_1x_0],[x_0x_1^{-1},x_0^{-2}x_1x_0^2]
\rangle$$ and infinitely as
$$
{\mathcal P} = \langle x_k, \ k \geq 0 | x_i^{-1}x_jx_i = x_{j+1}
\ \text{ if }i<j \rangle.$$
The infinite presentation provides a set of normal forms for elements of
$F$.  Namely, each $w \in F$ can be written
$x_{i_1}^{r_1}
x_{i_2}^{r_2}\ldots x_{i_k}^{r_k} x_{j_l}^{-s_l} \ldots
x_{j_2}^{-s_2} x_{j_1}^{-s_1} $ where $r_i, s_i >0$, and $i_1<i_2
\ldots < i_k$ and $j_1<j_2 \ldots < j_l$.   To obtain a unique normal
form for each element, we add the condition that when
 both $x_i$ and $x_i^{-1}$
occur, so does $x_{i+1}$ or  $x_{i+1}^{-1}$, as discussed by
Brown and Geoghegan
\cite{bg:thomp}.   We will always mean unique normal form
when we refer to a word $w$ in normal form.
We give a brief introduction to $F$ below; for a more detailed and
comprehensive description, we refer the reader to \cite{cfp} and
\cite{ct2}.

We interpret elements of $F$ geometrically as pairs of finite rooted binary
trees, each with the same number of exposed leaves, as described by Cannon, Floyd and
Parry in \cite{cfp}.  Let $T$ be a rooted binary tree.
An {\em exposed leaf} of $T$ ends in a vertex of valence $1$, and the exposed leaves are
numbered from left to right, beginning with $0$.
A node together with two leaves is called a {\em caret}.
A caret $C$ may have a {\em right child}, that is, a caret $C_R$ which is attached to the
right leaf of $C$.  We can similarly define the {\em left child} $C_L$ of the
caret $C$.   In a pair $(T_-,T_+)$ of rooted binary trees, the tree
$T_-$ is called the {\em negative tree} and $T_+$ the {\em positive tree}.
This terminology is explained further in \S \ref{sec:exponents} below.
The {\em left side} of
$T$ is the maximal path of left edges beginning at the root of
$T$. Similarly, we define the {\em right side} of $T$.

\begin{figure}
\includegraphics[width=.75in]{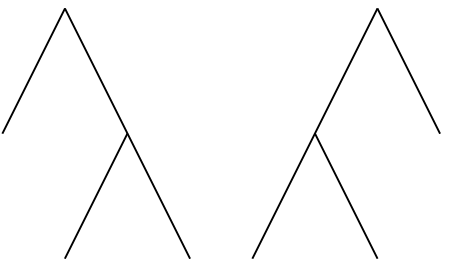} 
\caption{The tree pair diagrams for the generator $x_0$ of $F$, in the presentation $\F$.}
\label{fig:x0tree}

\includegraphics[width=.75in]{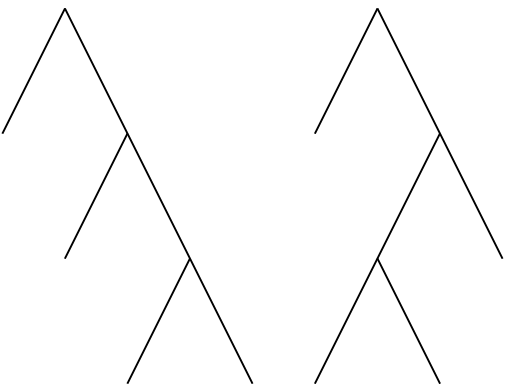} \\
\caption{ The tree pair diagram for the generator $x_1$ of $F$, in the presentation $\F$.}
\label{fig:x1tree}
\end{figure}

A pair of trees is {\em unreduced} if both $T_-$ and $T_+$
contain a caret with two exposed leaves numbered $m$ and $m+1$.
There are many tree pair diagrams representing the same element of
$F$ but each element has a unique reduced tree pair diagram representing it.
When we write $(T_-,T_+)$ to represent an element
of $F$, we are assuming that the tree pair is reduced. 
The tree pair diagrams for the generators $x_0$ and $x_1$ are given in figures \ref{fig:x0tree} and \ref{fig:x1tree}.

\subsection{Tree pair diagrams and the normal form}
\label{sec:exponents}

\begin{figure}
\label{fig:w}
\includegraphics[width=1.5in]{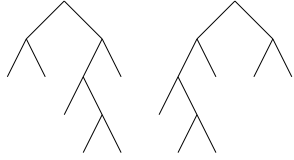} \\
\caption{ Tree pair diagram for the element $w=x_0^2 x_1 x_3^{-1} x_2^{-1} x_0^{-1}$.}
\label{fig:exptree}
\end{figure}

There is a bijective correspondence between the tree pair diagrams
described above and the normal form of an element.  The {\em leaf
exponent} of an exposed leaf numbered $n$ in $T_-$ or $T_+$ is
defined to be the length of the maximal path consisting entirely
of left edges from $n$ which does not reach the right side of the
tree, and is written $E(n)$. Note that $E(n)=0$ for an exposed
leaf labelled $n$ which is a right leaf of a caret, as there is no
path consisting entirely of left edges originating from $n$.
In figure \ref{fig:exptree}, the exponents of the left-hand tree $T_-$ are as follows: 
$E(0)=1, E(1)=0, E(2)=E(3)=1$ and $E(4)=E(5)=0$.

Once the exponents of the leaves in $T_-$ and $T_+$ have been computed, the
normal form of the element $w = (T_-,T_+)$ is easily obtained.  The positive
part of the normal form of $w$ is
$$
x_0^{E(0)} x_1^{E(1)} \cdots x_m^{E(m)}
$$
where $m$ is the number of exposed leaves in either tree, and the exponents are
obtained from the leaves of $T_+$.  The negative part of the normal form of $w$ is
similarly found to be
$$
x_m^{-E(m)} x_{m-1}^{-E(m-1)} \cdots x_0^{-E(0)}
$$
where the exponents are now computed from the leaves of $T_-$.
Note that many of the exponents in the normal form as given above
may be zero.

\subsection{Fordham's method of calculating word length}
\label{sec:Fordham}

Fordham \cite{blake} introduces a method of calculating, directly from the tree
pair diagram $w = (T_-,T_+)$, the word length $|w|$ of $w$ with respect to
the word metric arising from the finite presentation $\F$.  Throughout this paper we use $|w|$ to denote word length of a word $w$ in this word metric. In \cite{ct2} we use Fordham's method to estimate this word length in terms of the number of carets $N(w)$ in either tree of a tree pair diagram representing $w \in F$.

\begin{theorem}[\cite{ct2} theorem 3.1]
\label{thm:bound}
Let $w = (T_-,T_+)$ and $N(w)$ be the number of carets in $T_-$.  Then
$$N(w) - 2 \leq |w| \leq 4N(w) - 4.$$
\end{theorem}

Fordham's method of determining word length in $F$
 is based on a classification of the types of carets
found in a finite rooted binary tree $T$.  A rough classification
divides these carets into three types: {\em left carets}, which have a left
edge on the left side of $T$, {\em right carets}, which are not the root caret and
have a right edge on the right side of $T$, and {\em interior carets}.
Fordham specializes further by defining seven distinct types of carets, as follows.

\begin{enumerate}

\item
$L_0$.  The first caret on the left side of the tree, with
caret number $0$. Every tree has exactly one caret of type $L_0$.

\item
$\Ll$.  Any left caret other than the one numbered $0$.

\item
$I_0$. An interior caret which has no right child.

\item
$I_R$. An interior caret which has a right child.

\item
$R_I$. Any right caret numbered $k$ with the property that caret
$k+1$ is an interior caret.

\item
$\rni$. A right caret which is not an $R_I$ but for which there is a
higher numbered interior caret.

\item
$R_0$. A right caret with no higher-numbered interior carets.

\end{enumerate}

 The root caret is always considered to be a left caret and will
be of type $\Ll$ unless
it has no left children, in which case it would be the $L_0$ caret.

The carets in $T$ are numbered according to the
 infix ordering of nodes.  Caret $0$ is a left caret with an exposed left leaf
numbered $0$ in the leaf numbering.
According to the infix scheme, we number the left children of a caret before the caret
itself, and number the right children after numbering the caret.
Examples of trees whose carets are numbered in this way can be
found in \cite{ct}.

Fordham \cite{blake} proves that the word length $|w|$
of $w = (T_-,T_+)$ can be computed from knowing the caret types of
the carets in the two trees
 via the following process.  We number
the $k+1$ carets in each tree according to the infix method described above,
and for each $i$ with $0 \leq i \leq k$ we form the pair of caret
types consisting of the type of caret number $i$ in $T_-$ and the
type of caret number $i$ in $T_+$. The single caret of type $L_0$
in  $T_-$ will be paired with the single caret of type $L_0$ in
$T_+$, and for that pairing we assign a weight of 0. For all other
caret pairings, we assign weights according to the following
table.

\begin{center}

\begin{tabular}{|c|c|c|c|c|c|c|}

\hline
 & $R_0$ & $\rni$ & $R_I$ & $\Ll$ & $I_0$ & $I_R$ \\
 \hline

 $R_0$ & 0 & 2 & 2 & 1 & 1 & 3 \\ \hline
 $\rni$ & 2 & 2 & 2 & 1 & 1 & 3 \\ \hline
 $R_I$ & 2 & 2 & 2 & 1 & 3 & 3 \\ \hline
 $\Ll$ & 1 & 1 & 1 & 2 & 2 & 2 \\ \hline
 $I_0$ & 1 & 1 & 3 & 2 & 2 & 4 \\ \hline
$I_R$ & 3 & 3 & 3 & 2 & 4 & 4 \\ \hline
\end{tabular}

\end{center}

The main result of Fordham \cite{blake} is the following theorem.

\begin{theorem}[Fordham \cite{blake} 2.5.1]
\label{thm:blake} Given a word $w \in F$ described by the reduced tree
pair diagram $(T_-,T_+)$, the word length $|w|_{\F}$  is the sum of
the weights of the caret pairings in $(T_-,T_+)$.
\end{theorem}

\subsection{Group multiplication}
\label{sec:composition}

Group multiplication of $w = (T_-,T_+)$ and $v = (S_-,S_+)$ is
accomplished by creating temporary unreduced representatives
$(T'_-,T'_+)$ of $w$ and $(S'_-,S'_+)$ of $v$ in which $T'_+ =
S'_-$.  Then the product $wv$ is defined to be the tree pair
$(T'_-,S'_+)$ which may be unreduced. This method is used to
compute the distance between elements $w$ and $v$, namely, $d(w,v)
= |w^{-1} v|$.

\section{Quasi-isometric embeddings}
\label{sec:examples}

Let $X$ and $Y$ be metric spaces, and $K \geq 1$ and $C \geq 0$ be constants.
A {\em
$(K,C)$-quasi-isometric embedding} $f: X \rightarrow Y$ is a map
satisfying the following property:
$$ \frac{1}{K} d_X(x,y) - C \leq d_Y(f(x),f(y)) \leq Kd_X(x,y) + C$$
for all $x,y \in X$, where $d_X$ (resp. $d_Y$)
represents the metric in $X$ (resp. $Y$).  When considering a
quasi-isometric embedding between groups, we use the word metric
on each group induced by a particular set of generators.
A quasi-isometric embedding is a {\em quasi-isometry} if there
is a constant $C'$ so that $Nbhd_{C'}(f(X)) = Y$.

The definition of a quasi-isometry does not include any algebraic
requirements, in fact, a quasi-isometry need not be a continuous
map.  In addition, a quasi-isometry is not required to preserve
any algebraic structure at all.  Below we give several examples of
quasi-isometric embeddings of $F^m  \times \Z^n$ into $F$, for nonnegative integers $m$ and $n$, which, though homomorphisms,
are algebraically cumbersome.  
However, they are easily described geometrically by the way they alter the tree pair diagram of an element $w \in F$.  
We call these {\em geometric}
quasi-isometric embeddings.

Let $G$ be a finitely generated group and $H$ a finitely generated
subgroup of $G$.  One can ask if the inclusion of $H$ into $G$ is
a quasi-isometric embedding.  The inclusion will be a
quasi-isometric embedding if and only if  the distortion function of the group
is bounded. Define the {\em distortion function} $h(n)$ of $H$ in
$G$ by
$$
h(n) = \frac{1}{n} \max \{ |x|_H  \ | x \in H \text{ and } |x|_G \leq n \}.
$$
We will use the distortion function to show that the examples we
give below are quasi-isometric embeddings.

\subsection{The shift map $\phi$}
\label{sec:phi}

We begin with an example of a quasi-isometric embedding $F \rightarrow F$
which is easily understood either algebraically or geometrically.
There is a {\em shift map} $\phi: F \rightarrow F$ defined on Thompson's
group in the infinite presentation ${\mathcal P}$ which increases the index of each
generator by $1$.  For example, if $w = x_3^2 x_5 x_{13} x_{10}^{-1}
x_9^{-4}$, then $\phi(w) = x_4^2 x_6 x_{14} x_{11}^{-1} x_{10}^{-4}$.
Thus if $w$ is in normal form, so is $\phi(w)$.
We show that any power of $\phi$ is a quasi-isometric embedding.

\begin{figure}
\includegraphics[width=1.5in]{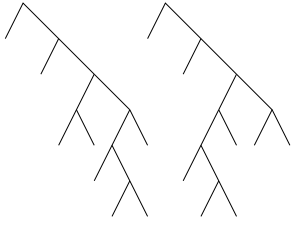} \\
\caption{Tree pair diagram for $\phi^2(w)=x_2^2 x_3 x_5^{-1} x_4^{-1} x_2^{-1}$, where the word $w=x_0^2 x_1 x_3^{-1} x_2^{-1} x_0^{-1}$ is depicted in figure 3.}
\label{fig:phi2w}
\end{figure}

\begin{theorem}
\label{thm:phi}
Any integral power $\phi^n:F \rightarrow F$ of the shift map $\phi: F \rightarrow F$
is a quasi-isometric embedding.
\end{theorem}

\begin{proof}
We first note that $\phi^n(F)$ is isomorphic to the subgroup of
$F$ generated by $x_n$ and $x_{n+1}$.  We will therefore use
distortion to show that $\phi^n$ is a quasi-isometric embedding.

Let $w=(T_-,T_+)$, and $\phi^n(w) = (S_-,S_+)$.  We first describe
the relationship between the trees $T_{\pm}$ and $S_{\pm}$.  The
smallest possible index of a generator in $\phi^n(w)$ is $n$.
Thus, the leaf exponents of all leaves numbered $0$ through $n-1$
in $S_{\pm}$ must be $0$.  The trees $S_{\pm}$ are now easily
deduced from $T_{\pm}$. The tree $S_-$ must have an empty left
subtree, $n-1$ right carets with exposed left carets beginning at
the root, and then the tree $T_-$.  Similarly, $S_+$ is the tree
having an empty left subtree, $n-1$ right carets with no left
subtrees to the right of the root, and then the tree $T_+$.  
Figure \ref{fig:phi2w} gives an example of he tree pair diagram of $\phi^2(w)$, using the word $w$ depicted in figure \ref{fig:exptree}.  

We now show that $\phi^n$ is a quasi-isometric embedding. Recall
from theorem \ref{thm:bound} that 
$$N(w) - 2 \leq |w| \leq 4N(w) -4,$$ 
where $N(w)$ is the number of carets in $T_-$. Since
$N(\phi^n(w)) = N(w) + n$, we obtain the same inequality for
$|\phi^n(w)|$:
$$
N(w) + (n-2) \leq |\phi^n(w)| \leq 4(N(w) + n) - 4 = 4N(w) +
4(n-1).
$$
Combining the two above inequalities yields the desired bound:
$$
\frac{1}{4}  |w| + (n+2) \leq |\phi^n(w)| \leq 4|w| + (4n-2).
$$
Since $n$ is a constant, it follows that the distortion function is
bounded and thus $\phi^n$ is a quasi-isometric embedding.
\end{proof}

\subsection{The reverse shift map $\psi$}
\label{sec:psi}
\begin{figure}
\includegraphics[width=1.5in]{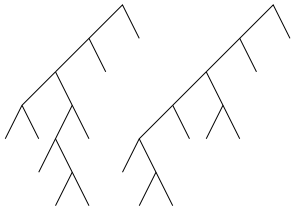} \\
\caption{Tree pair diagram for $\psi^2(w)=x_0^4 x_1 x_4 x_3^{-1} x_2^{-2} x_0^{-3}$, where the word $w=x_0^2 x_1 x_3^{-1} x_2^{-1} x_0^{-1}$ is depicted in figure 3.}
\label{fig:psi2w}
\end{figure}

Given the simple geometric representation of the shift map $\phi$, it is
natural to define a ``reverse shift map'', which we denote $\psi$, in which the original trees become left subtrees of the root caret rather than right subtrees.  See, for example, figure \ref{fig:psi2w}, which depicts the tree pair diagram for $\psi^2(w)$.  Geometrically, it is completely natural to
consider $\psi$ as well as $\phi$.  However, in the literature,
$\psi$ is rarely mentioned, likely because it is algebraically cumbersome,
unlike $\phi$, in the following way.
If $w \in F$ is written in normal form, then the normal form of
$\phi(w)$ is easily determined, whereas for $\psi$ the situation
is quite different. Let $w = (T_-,T_+)$ and $\psi(w) = (S_-,S_+)$.
Then leaves with exponent $0$ which are the exposed left leaves of
right carets in $T_-$ or $T_+$ will have leaf exponent $1$ in
$S_-$ and $S_+$, respectively, and thus cause new generators to
appear in the normal form of $\psi(w)$. It is difficult to predict
solely from the normal form of $w$ which additional generators
will appear in the normal form of $\psi(w)$. For example if $w =
x_0^2 x_1 x_5^2 x_4^{-1} x_3^{-1} x_1^{-1} x_0^{-1}$ then $\psi(w)
= x_0^3 x_1 x_4 x_5^3 x_7^{-1} x_6^{-1} x_4^{-1} x_3^{-2} x_1^{-1}
x_0^{-1}$.

Unsurprisingly from the geometric point of view, we
obtain a theorem for $\psi$ analogous to theorem \ref{thm:phi} for $\phi$.

\begin{theorem}
\label{thm:psi}
Any integral power $\psi^n:F \rightarrow F$ of the reverse shift map
$\psi: F \rightarrow F$ defined above is a quasi-isometric embedding.
\end{theorem}

\begin{proof}
To apply the method of proof used above in theorem \ref{thm:phi},
we give the generators of a subgroup of $F$ isomorphic to
$\psi^n(F)$.  Once we do this, it is clear from the tree pair diagrams
 for $w$ and $\psi^n(w)$ defined above that $S_{\pm}$ has
$n$ more carets than $T_{\pm}$, and we obtain inequalities for
$|\psi^n(w)|$ identical to those for $|\phi^n(w)|$ given in the
proof of theorem \ref{thm:phi}.  Thus $\psi^n$ will be a quasi-isometric
embedding of $F$ into $F$.

It is easily checked that $\psi$ is a group homomorphism, and thus
$\psi^n(F)$ is generated by $\psi^n(x_0)$ and $\psi^n(x_1)$. After
forming the tree pair diagrams representing these generators,
using the exponent notation of \S \ref{sec:exponents}, we see that
$\psi^n(\x) = \x^{n} x_1 \x ^{-(n+1)}$ and $\psi(x_1) = \x^n x_1^2
x_2^{-1} x_1^{-1} \x^{-(n+1)}$
which gives  an algebraic
description of $\psi$. Similarly, we can compute $$\psi^n(x_k) =
x_0^n x_1 x_2 \cdots x_{k-1} x_k^2 x_{k+1}^{-1} x_k^{-1} \cdots
x_2^{-1} x_1^{-1} x_0^{-n}.$$ Substituting this formula into the
normal form of $w \in F$, we obtain an expression for $\psi^n(w)$
in the infinite generating set which can be simplified to obtain
the normal form for $\psi^n(w)$.
\end{proof}

\subsection{Clone subgroups of $F$}
\label{sec:clone}

The maps $\phi^n$ and $\psi^n$ are shown to be
quasi-isometric embeddings because an element and its image are
represented by tree pair diagrams which differ by a finite number of carets. 
This idea can be generalized as
follows.   Let $p$ denote the composition 
$$p=f_1 \circ f_2 \circ \cdots \circ f_n$$ where each $f_i$ is either $\phi$ or $\psi$ and we compose from left to right, so that $f_1 \circ f_2(F) = f_1(f_2(F))$.  A {\em clone subgroup} is the image of $F$ under $p$, that is, $p(F)$.  The simplest examples of clone subgroups are $\phi^n(F)$ and $\psi^n(F)$.  Clone subgroups can be understood in a number of ways.

\begin{figure}
\includegraphics[width=1in]{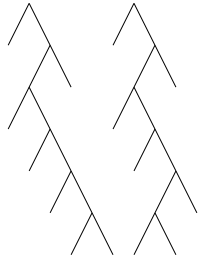} \\
\caption{Let $p=\phi \circ \psi \circ \phi^2$. This tree pair represents the element  $p(x_0)$ in the corresponding clone subgroup $p(F)=C_{1011}$.}
\label{fig:clonex0tree}
\end{figure}

\begin{figure}
\includegraphics[width=1in]{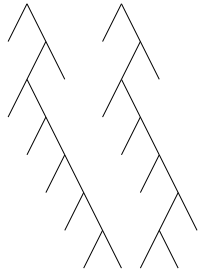} \\
\caption{ Let $p=\phi \circ \psi \circ \phi^2$. This tree pair represents the element  $p(x_1)$ in the clone subgroup $p(F)=C_{1011}$.}
\label{fig:clonex1tree}
\end{figure}

Let $p$ be as above, and consider the  clone subgroup $p(F)$. We can describe this subgroup using a binary address. A node in a binary tree is given an address using the following inductive method.  
The root node has the empty label.  Given a node with label $s$,  the left child of the node
is labelled $s0$ and  right child of the node labelled $s1$.   For example, the right child of the right child of the left child of the right child of the root has address $1011$.

Given $p$ as above, and $w=(T_-,T_+) \in F$, let $p(w) = (S_-,S_+)$. From the definitions of $\phi, \ \psi$ and $p$, we know that the tree $T_-$ will be a subtree of $S_-$, and $T_+$ will be a subtree of $S_+$. We consider the address $s= \epsilon_1 \epsilon_2 \ldots \epsilon_n$ of the root caret of $T_-$ as a subtree of $S_-$ and use $s$ to give an ``address" for the clone subgroup $p(F)$.  Namely, we can uniquely refer to $p(F)$ as $C_s$.  
For example, the subgroup $C_{1011}$ is the
image $\phi(\psi (\phi^2(F)))$.

We can also describe clone subgroups by their representation as piecewise linear homeomorphisms
of the unit interval.  Each dyadic subinterval of the form $[\frac{i}{2^n},\frac{i+1}{2^n}]$ is equivalent to the standard unit interval $[0,1]$ by an affine map with dyadic coefficients.  For a fixed clone subgroup $C_s$ where $s$ has length $n \geq 1$, there is  a proper subinterval
$I' $ of $ [0,1]$ which contains the $x$- and $y$- coordinates of all the breakpoints of elements of $C_s$, 
and such that all elements of $C_s$ are the identity map outside of the interval $I'$.
The length of $I'$ is $2^{-n}$ and the endpoints of $I'$are $\frac{i}{2^n}$ and $\frac{i+1}{2^n}$.   The endpoints can be computed easily from the address $s$. 
Furthermore, any element
of $F$ whose breakpoints all lie in the interval $[\frac{i}{2^n}, \frac{i+1}{2^n}]$ will be an element of $C_s$.
For the example subgroup $C_{1011}$, the  dyadic interval containing the
all breakpoints of its elements is $[\frac{11}{16},\frac34]$.

It is easy to see that any clone subgroup $C_s=p(F)$ is isomorphic to $F$.  
If $w = (T_-,T_+)$ and $p(w) = (S_-,S_+)$, we know by definition that the tree pair $(T_-,T_+)$ is reduced.  From the definitions of $\phi$ and $\psi$, and thus $p$, no additional reduction can occur when the tree pair $(S_-,S_+)$ is formed.  
Thus each element of $F$ produces a unique element of $C_s$, and it is clear that $p(x_0)$ and $p(x_1)$ must generate $C_s$.  

The following theorem is a consequence of theorems \ref{thm:phi} and \ref{thm:psi}.
\begin{theorem}
Any clone subgroup of $F$ is quasi-isometrically embedded.
\end{theorem}

\begin{proof} The clone subgroup $C_s$ is determined by $p=f_1 \circ f_2 \circ \cdots \circ f_n$ where each $f_i$ is either the map $\phi$ or $\psi$.  
We know from theorems \ref{thm:phi} and \ref{thm:psi} that $\phi$ and $\psi$ are $(4,2)$-quasi-isometric embeddings of $F$ into $F$. Geometrically, we see that if $w=(T_-,T_+)$ and
$p(w)=(S_-,S_+)$, the trees $T_\pm$ and $S_\pm$ differ by $n$ carets. Thus,
we again get that $\frac{1}{4}  |w| + (n+2) \leq |p(w)| \leq 4|w| + (4n-2)$ for all $w$ in $F$, which shows that $p$ is also a quasi-isometric embedding, with constants $(4,4n-2)$.  Thus $C_s$ is quasi-isometrically embedded.  
\end{proof}

\subsection{Quasi-isometrically embedded products}
We now show that $F$ contains a family of quasi-isometrically embedded subgroups of the form $F^m \times \Z^n$ for $n,m \geq 0$.  These examples include those of Burillo \cite{burillo}, who presents a family of subgroups of $F$ of the form $F \times \Z^n$ which are quasi-isometrically embedded.     Guba and Sapir \cite{diag,diag2}, using the approach of diagram groups,  also show that
$F^m \times \Z^n$ occur as quasi-isometrically embedded subgroups of $F$ and furthermore
show that all abelian subgroups and centralizers of elements are quasi-isometrically
embedded.  We present an alternate geometric proof that  some concrete geometric
realizations of the subgroups $F^m \times \Z^n$ are quasi-isometrically embedded,
using the shift maps described above. 

\begin{theorem}
\label{thm:products}
For each integral pair $m,n \geq 0$, Thompson's group $F$ contains an infinite family of quasi-isometrically embedded subgroups of the form $F^m \times \Z^n$.
\end{theorem} 

As mentioned above,  theorem \ref{thm:products} generalizes the following theorem and corollary of Burillo \cite{burillo}.

\begin{figure}
\includegraphics[width=3.5in]{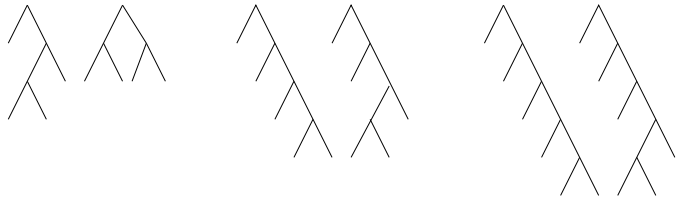} \\
\caption{ The generators $x_0x_1^{-1}, \ x_2$ and $x_3$, respectively, of the quasi-isometrically embedded $F \times \Z$ described in theorem \ref{thm:burillo}.}
\label{fig:fxz}
\end{figure}

\begin{theorem}[Burillo \cite{burillo}, theorem 3]
\label{thm:burillo}
The map $\Phi: F \times \Z \rightarrow F$ defined by $\Phi(x,t) = (x_0x_1^{-1})^t
\phi^2(x)$ where $x \in F$ and $t \in \Z$, is a quasi-isometric embedding.
\end{theorem}

\begin{corollary}[Burillo \cite{burillo}, corollary 4]
\label{cor:burillocor}
The subgroup generated by the elements $$x_0x_1^{-1}, \
x_2x_3^{-1}, \cdots x_{2n-2}x_{2n-1}^{-1}, \ x_{2n}, \
x_{2n+1}$$ is isomorphic to $F \times \Z^n$ and is quasi-isometrically embedded.
\end{corollary}

The first $n$ elements in corollary \ref{cor:burillocor}
generate a subgroup isomorphic to $\Z^n$, and the final two generate a subgroup isomorphic
to $F$.

Burillo proves these results by bounding the distortion function
using a technique he develops to estimate the length of $w$ in
the word metric corresponding to the generating set of the presentation $\F$ as a
function of the normal form of $w$.  We present an alternate geometric proof of this theorem and corollary before presenting the proof of theorem \ref{thm:products}.   We note that
Guba and Sapir \cite{diag2} show that many subgroups of the form $F^m \times \Z^n$
are quasi-isometrically embedded using methods for diagram groups.

Consider the quasi-isometric embedding $\Phi$ of theorem
\ref{thm:burillo}. We identify geometrically the subgroup $F
\times \Z$ of $F$ generated by the elements $x_0 x_1^{-1}, \ x_2$
and $x_3$. Drawing the tree pairs representing $x_2$ and
$x_3$, we see that they generate a clone subgroup $C_{11}$ of
$F$.  The tree pair diagrams representing these generators are given in figure \ref{fig:fxz}.  

An element of $C_{11}$, when viewed as an element of the larger group
$F$, has no instances of $\x^{\pm 1}$ or $x_1^{\pm 1}$ in its
normal form, and therefore lies in the subgroup $\phi^2(F)$.  We
also know that $\phi^2(F)$ is generated by $\phi^2(\x) = x_2$ and
$\phi^2(x_1) = x_3$.

It is easily seen that $(x_0x_1^{-1})^k = x_0^k x_k^{-1} x_{k-1}^{-1} \cdots
x_2^{-1} x_1^{-1}$, and we draw the tree pair diagram $(S_-,S_+)$ representing
this element as follows.  The form of the tree $S_+$ is clear, and in determining
$S_-$, we see that the left subtree of the root caret is empty, and the
right subtree of the root caret has as string of $k-1$ carets of type $I_R$, a single
caret of type $I_0$ and a
single caret of type $R_0$.

We are really interested in the fact that the right subtree of the
right child of the root caret in $T_{\pm}$ is empty.  This can be
easily determined by applying the following lemma from \cite{ct2}
to the substring  $x_k^{-1} x_{k-1}^{-1} \cdots x_2^{-1} x_1^{-1}$
of $(x_0x_1^{-1})^k $, since this word  determines the right subtree of the
root of $S_-$.

\begin{lemma}[\cite{ct2}, lemma 5.1]
\label{lemma:rightside}
Let the element $w = (T_-,T_+)$ have normal form
$$x_0^{r_0} x_{i_1}^{r_1} x_{i_2}^{r_2}\ldots x_{i_k}^{r_k} x_{j_l}^{-s_l}
\ldots x_{j_2}^{-s_2} x_{j_1}^{-s_1} x_0^{-s_0}.$$

\begin{enumerate}

\item
If $\displaystyle i_k < \Sigma_{m=0}^{k-1} r_m$  then the right subtree of
the root caret of $T_+$ is empty.

\item
If $\displaystyle j_l < \Sigma_{m=0}^{l-1} s_m$  then the right subtree of
the root caret of $T_-$ is empty.
\end{enumerate}
\end{lemma}

Since $k< \Sigma_{i=1}^{k-1} i$, we see that the right subtree of the
right child of the root caret of $T_-$ is empty, as is the
corresponding subtree in $T_+$. This fact guarantees that an
element in $C_{11} = \langle x_2, \ x_3 \rangle$ and $(x_0x_1^{-1})^k$
commute, as is easily seen from the definition of group multiplication of pairs of trees,
or simply by noting that subintervals where these elements are different from
the identity map are disjoint.

We apply the method of counting carets used above in \S \ref{sec:phi} and
\ref{sec:psi} to give an alternate proof of theorem \ref{thm:burillo}.
We first note that $C_{11}$ is the image $\phi^2(F)$
It is easily seen that $N(\phi^2(w)) = N(w) + 2$ and $N((x_0x_1^{-1})^t)
=t+2$.  The two parts of $w$ are in different subtrees, excepting
an overlap of two carets at the root and the right child of the root.
So when we count carets, we get $N(\Phi(w,t)) = N(w) + t + 2$.
Thus
$$
|\Phi(w,t)|_{\F} \leq 4 N(\Phi(w,t))-4 \leq 4 (N(w)+t+2) -4 \leq 
4|w| + 4t + 12 \leq 4(|(w,t)|_{F \times \Z} ) + 12
$$
and
$$
|\Phi(w,t)|_{\F} \geq N(\Phi(w,t)) -2 \geq N(w) +t \geq \frac{1}{4} |w| + t + 4 \geq
\frac{1}{4} (|(w,t)|_{F \times \Z} ) + 4.
$$
Thus the distortion function is bounded and the map $\Phi$ is a
quasi-isometric embedding.

To extend this proof to a geometric interpretation of corollary \ref{cor:burillocor}, notice that the first $n-1$ terms listed in the corollary, $x_0x_1^{-1}, \cdots ,x_{2n-2}x_{2n-1}^{-1}$, commute with each other.  This is easily seen from the tree pair diagrams representing these elements or simply by noting that these elements lie in clone subgroups that intersect only in the identity.
The last two generators of the corollary, $x_{2n}$ and $x_{2n+1}^{-1}$, generate a clone subgroup $C_{11\ldots11}$, whose address is $2n$ one's.

\medskip
\noindent
{\it Proof of theorem \ref{thm:products}.}
Fix  a pair of integral $m,n \geq 0$.  We now produce generators for  subgroups of the form  $F^m \times \Z^n$ of $F$ and show that these subgroups are quasi-isometrically embedded.  Each copy of $F$ comprising the $F^m$ will be given by a distinct clone subgroup.  We can let $s_1,s_2, \cdots s_m, s_{m+1}$ be binary addresses of nodes in a binary rooted tree,  subject to the condition that pairwise, no $s_i$ is a prefix of $s_j$.  For a fixed $m$ there are
infinitely many choices of prefix sets satisfying this requirement.
Then for each $i$ with $1 \leq i \leq m$, we have that $C_{s_i}$ is a clone subgroup of $F$, and by construction, no two of these clone subgroups intersect in  non-identity element.  Thus the product $C_{s_1} \times \cdots \times C_{s_m}$ is isomorphic to $F^m$.  The final factor clone subgroup $C_{s_{m+1}}$ will contain the abelian factors.

Let $p_{m+1}$ be the product of maps $\phi$ and $\psi$ corresponding to the binary address $s_{m+1}$, and thus $p_{m+1}(F) = C_{s_{m+1}}$.  Using the elements of corollary \ref{cor:burillocor}, we see that $$p_{m+1}(x_0x_1^{-1}), \ p_{m+1}(x_2x_3^{-1}), \cdots , p_{m+1}(x_{2n-2}x_{2n-1}^{-1})$$ will generate a subgroup of $F$ isomorphic to $\Z^n$.  

We now claim that the map $\Psi: F^m \times \Z^n \rightarrow F$ given by 
$$\Psi(w_1,...,w_m,t_1,...,t_n)=p_{s_1}(w_1) \cdots p_{s_m}(w_m)(p_{s_{m+1}}(x_0x_1^{-1}))^{t_1} \cdots (p_{s_{m+1}}(x_{2n-2}x_{2n-1}^{-1}))^{t_n}
$$
is a quasi-isometric embedding.  We again apply the method of counting carets to prove the claim.  

First we note that all the factors have been chosen to lie in clone subgroups which intersect
only in the identity, so they commute.  Since each factor lies in a clone subgroup $C_{s_i}$  there are constants $k_i$ so that $N(p_{s_i}(w_i)) = N(w_i) + k_i$, where $p_{s_i}$ is product of $k_i$ maps.  
Similarly, from the tree pair diagrams, it is easy to see that $N((x_{2i}x_{2i+1}^{-1})^t) = N(x_{2i}x_{2i+1}^{-1})+t$.  
Thus $$N(p_{s_{m+1}}(x_{2i}x_{2i+1}^{-1})^{t_i}) = N(p_{s_{m+1}}(x_{2i}x_{2i+1}^{-1}))+t_i = N(x_{2i}x_{2i+1}^{-1})+k_{m+1} + t_i,$$

 Combining the above inequalities with theorem \ref{thm:bound}, we obtain the following upper bound.
$$
|\Psi(w_1,w_2, \cdots ,w_m,t_1,t_2, \cdots ,t_n)|_{\F} 
$$
$$\leq 4 \Big( \Sigma_{i=1}^m (|w_i|_{\F} + k_i) + \Sigma_{j=1}^n (|x_{2j}x_{2j+1}^{-1}|_{\F} + t_i+k_{m+1})\Big)  -4$$
$$
= 4\Sigma_{i=1}^m |w_i|_{\F} + 4\Sigma_{j=1}^n t_j +C' 
$$
$$\leq 4|(w_1,w_2, \cdots ,w_m,t_1,t_2, \cdots ,t_n)|_{F^m \times \Z^n} +C'$$
where $C' = 4(\Sigma k_i + 4\Sigma |x_{2j}x_{2j+1}^{-1}|_{\F} + nk_{m+1} -1) $.
The lower bound 
$$
|\Psi(w_1,w_2, \cdots ,w_m,t_1,t_2, \cdots ,t_n)|_{\F} \geq |(w_1,w_2, \cdots ,w_m,t_1,t_2, \cdots ,t_n)|_{F^m \times \Z^n} -C''
$$
is obtained analogously, for the appropriate constant $C''$ depending on $s_1, \cdots ,s_{m+1}$.  

Thus the map $\Psi$ is a quasi-isometric embedding.  Since there are infinitely many choices for $s_1, \cdots ,s_{m+1}$, we have produced an infinite family of quasi-isometrically embedded subgroups of $F$ of the form $F^m \times \Z^n$ for any integral pair $(m,n)$ with $m,n\geq0$.  
\qed

\bigskip

\begin{small}
\noindent Sean Cleary \\
Department of Mathematics \\
City College of New York \\
City University of New York \\
New York, NY 10031 \\
E-mail: {\tt cleary@sci.ccny.cuny.edu}

\medskip

\noindent
Jennifer Taback\\
Department of Mathematics and Statistics\\
University at Albany\\
Albany, NY 12222\\
E-mail: {\tt jtaback@math.albany.edu}
\end{small}

\end{document}